\documentclass[12pt]{article}%
\usepackage{amsmath}
\usepackage{amsfonts}
\usepackage{amssymb}
\usepackage{graphicx}%
\newtheorem{theorem}{Theorem}

\newtheorem{corollary}[theorem]{Corollary}

\newtheorem{definition}[theorem]{Definition}
\newtheorem{example}[theorem]{Example}

\newtheorem{lemma}[theorem]{Lemma}

\newtheorem{proposition}[theorem]{Proposition}

\newenvironment{proof}[1][Proof]{\noindent\textbf{#1.} }{\ \rule{0.5em}{0.5em}}
\begin{document}

\title{The Classification of All Singular Nonsymmetric Macdonald Polynomials}
\author{Charles F. Dunkl\thanks{Email: cfd5z@virginia.edu}\\Dept. of Mathematics\\University of Virginia\\Charlottesville VA 22904-4137}
\maketitle

\begin{abstract}
The affine Hecke algebra of type $A$ has two parameters $\left(  q,t\right)  $
and acts on polynomials in $N$ variables. There are two important pairwise
commuting sets of elements in the algebra: the Cherednik operators and the
Jucys-Murphy elements whose simultaneous eigenfunctions are the nonsymmetric
Macdonald polynomials, and basis vectors of irreducible modules of the Hecke
algebra, respectively. For certain parameter values it is possible for special
polynomials to be simultaneous eigenfunctions with equal corresponding
eigenvalues of both sets of operators. These are called singular polynomials.
The possible parameter values are of the form $q^{m}=t^{-n}$ with $2\leq n\leq
N.$ For a fixed parameter the singular polynomials span an irreducible module
of the Hecke algebra. Colmenarejo and the author (SIGMA 16 (2020), 010) showed
that there exist singular polynomials for each of these parameter values, they
coincide with specializations of nonsymmetric Macdonald polynomials, and the
isotype (a partition of $N$) of the Hecke algebra module is $\left(
dn-1,n-1,\ldots,n-1,r\right)  $ for some $d\geq1$. In the present paper it is
shown that there are no other singular polynomials.

\end{abstract}

\section{Introduction}

Many structures arise from the action of the symmetric group on polynomials in
$N$ variables. Among them are the Hecke algebra and the affine Hecke algebra
of type $A$. This paper concerns polynomials with noteworthy properties with
respect to these algebras. The symmetric group $\mathcal{S}_{N}$ is generated
by the simple reflections $s_{i},1\leq i<N$,  where
\[
xs_{i}:=\left(  x_{1},\ldots,\overset{i}{x_{i+1}},\overset{i+1}{x_{i}}%
,\ldots,x_{N}\right)  ;
\]
they satisfy the braid relations $s_{i}s_{i+1}s_{i}=s_{i+1}s_{i}s_{i+1}$ and
$s_{i}s_{j}=s_{j}s_{i}$ for $\left\vert i-j\right\vert \geq2$. Let $q,t$ be
parameters satisfying $t^{n}\neq1$ for $2\leq n\leq N$ and $q,t\neq0$. Define
$\mathcal{P=}\mathbb{K}\left[  x_{1},\ldots,x_{N}\right]  $ where $\mathbb{K}$
is a field containing $\mathbb{Q}\left(  q,t\right)  $. The Hecke algebra
$\mathcal{H}_{N}\left(  t\right)  $ is generated by Demazure operators (with
$p\in\mathcal{P}$ and $1\leq i<N$)%
\[
T_{i}p\left(  x\right)  :=\left(  1-t\right)  x_{i+1}\frac{p\left(  x\right)
-p\left(  xs_{i}\right)  }{x_{i}-x_{i+1}}+tp\left(  xs_{i}\right)  ;
\]
they satisfy the same braid relations $T_{i}T_{i+1}T_{i}=T_{i+1}T_{i}T_{i+1}$
and $T_{i}T_{j}=T_{j}T_{i}$ for $\left\vert i-j\right\vert \geq2$, as well as
the quadratic relations $\left(  T_{i}-t\right)  \left(  T_{i}+1\right)  =0$.
The affine Hecke algebra $\mathcal{H}_{N}\left(  t;q\right)  $ is obtained by
adjoining the $q$-shift%
\[
wp\left(  x\right)  :=p\left(  qx_{N},x_{1},x_{2},\ldots,x_{N-1}\right)
\]
and defining%
\begin{align*}
T_{0}p\left(  x\right)   &  :=wT_{1}w^{-1}p\left(  x\right)  =\left(
1-t\right)  x_{1}\frac{p\left(  x\right)  -p\left(  xs_{0}\right)  }%
{qx_{N}-x_{1}}+tp\left(  xs_{0}\right)  ,\\
xs_{0} &  :=\left(  qx_{N},x_{2},\ldots,x_{N-1},x_{1}/q\right)  .
\end{align*}
Then $wT_{i+1}=T_{i}w$ where the indices are taken $\operatorname{mod}N$.
(That is, $w^{2}T_{1}=T_{N-1}w^{2}$.) The quadratic relations imply
$T_{i}^{-1}=t^{-1}\left(  T_{i}+\left(  1-t\right)  \right)  $. There are two
commutative families of operators in $\mathcal{H}_{N}\left(  t;q\right)  $
(each indexed $1\leq i\leq N$): the Cherednik operators (see \cite{C1995})%
\[
\xi_{i}:=t^{i-1}T_{i}T_{i+1}\cdots T_{N-1}wT_{1}^{-1}T_{2}^{-1}\cdots
T_{i-1}^{-1}%
\]
and the Jucys-Murphy operators
\[
\omega_{i}=t^{i-N}T_{i}T_{i+1}\ldots\cdots T_{N-1}T_{N-1}T_{N-2}\cdots T_{i}.
\]
Note that $\xi_{i}=t^{-1}T_{i}\xi_{i+1}T_{i}$ and $\omega_{i}=t^{-1}%
T_{i}\omega_{i+1}T_{i}$ for $i<N$. The simultaneous eigenfunctions of the
Cherednik operators are the nonsymmetric Macdonald polynomials and the
simultaneous eigenvectors of the Jucys-Murphy operators span irreducible
representations of $\mathcal{H}_{N}\left(  t\right)  $. Our concern is to
determine all polynomials which are simultaneous eigenfunctions of both sets
of operators, more specifically, when $q,t$ satisfy a relation of the form
$q^{m}t^{n}=1$ to determine the homogeneous polynomials $p$ such that $\xi
_{i}p=\omega_{i}p$ for all $i$. These are called \textit{singular polynomials}
with singular parameter $q^{m}=t^{-n}$. In a previous paper \cite{CD2020}
Colmenarejo and the author found a large class of such polynomials associated
with tableaux of quasi-staircase shape. In this paper we will show that there
are no other occurrences.

Affine Hecke algebras were used by Kirillov and Noumi \cite{KN1998} to derive
important results about the coefficients of Macdonald polynomials. Mimachi and
Noumi \cite{MN1998} found double sums for reproducing kernels for series in
nonsymmetric Macdonald polynomials. The paper \cite{BF1997} by Baker and
Forrester is a source of some background for the present paper.

In Section \ref{prelim} we collect the needed definitions and results about
the Hecke algebra action on polynomials, Cherednik operators, nonsymmetric
Macdonald polynomials, and the representation theory of Hecke algebra of type
$A$. The definition of singular polynomials and its consequences, that is,
necessary conditions, are presented in Section \ref{necessary}. This section
also explains the known existence theorem. Section \ref{restrict} concerns the
method of restriction to produce singular polynomials with a smaller number of
variables and this leads into Section \ref{conclud} where our main
nonexistence theorem is proved.

\section{\label{prelim}Preliminary Results}

In this section we present background information and computational results
dealing with $\mathcal{H}_{N}\left(  t\right)  $ and the action on polynomials.

\begin{lemma}
\label{comrel1}If $j>i+1$ or $j<i$ then $T_{i}\omega_{j}=\omega_{j}T_{i}$, and
$T_{i}\omega_{i}=\left(  t-1\right)  \omega_{i}+\omega_{i+1}T_{i}$,
$T_{i}\omega_{i+1}=\omega_{i}T_{i}-\left(  t-1\right)  \omega_{i}$.
\end{lemma}

\begin{proof}
If $j>i+1$ then $T_{i}$ commutes with each factor of $\omega_{i}$. Suppose
$j=i-1$ then by the braid relations
\begin{align*}
T_{i}\omega_{i-1}  &  =t^{i-1-N}T_{i}T_{i-1}T_{i}T_{i+1}\cdots T_{i}%
T_{i-1}=t^{i-1-N}T_{i-1}T_{i}T_{i-1}T_{i+1}\cdots T_{i}T_{i-1}\\
&  =t^{i-1-N}T_{i-1}T_{i+1}T_{i+2}\cdots T_{i-1}T_{i}T_{i-1}\\
&  =t^{i-1-N}T_{i-1}T_{i}T_{i+1}\cdots T_{i}T_{i-1}T_{i}=\omega_{i-1}T_{i}.
\end{align*}
Suppose $j<i-1$ then $\omega_{j}=t^{j-i+1}T_{j}T_{j+1}\cdots T_{i-2}%
\omega_{i-1}T_{i-2}\cdots T_{j}$ and $T_{i}$ commutes with each factor in this
product. If $j=i$ then
\begin{align*}
T_{i}\omega_{i}  &  =t^{-1}T_{i}^{2}\omega_{i+1}T_{i}=t^{-1}\left\{  \left(
t-1\right)  T_{i}+t\right\}  \omega_{i+1}T_{i}\\
&  =\left(  t-1\right)  \omega_{i}+\omega_{i+1}T_{i},
\end{align*}
and similarly $\omega_{i}T_{i}=t^{-1}T_{i}\omega_{i+1}T_{i}^{2}=\left(
t-1\right)  \omega_{i}+T_{i}\omega_{i+1}$
\end{proof}

\begin{lemma}
\label{comrel2}If $j>i+1$ or $j<i$ then $T_{i}\xi_{j}=\xi_{j}T_{i}$, and
$T_{i}\xi_{i}=\left(  t-1\right)  \xi_{i}+\xi_{i+1}T_{i},T_{i}\xi_{i+1}%
=\xi_{i}T_{i}-\left(  t-1\right)  \xi_{i}$.
\end{lemma}

\begin{proof}
Recall $wT_{i+1}=T_{i}w,~w^{2}T_{1}=T_{N-1}w^{2}$. Suppose $j=i-1$ then%
\begin{align*}
T_{i}\xi_{i-1}  &  =t^{i-1-N}T_{i}T_{i-1}T_{i}T_{i+1}\cdots T_{N-1}wT_{1}%
^{-1}\cdots T_{i-2}^{-1}\\
&  =t^{i-1-N}T_{i-1}T_{i}T_{i-1}T_{i+1}\cdots T_{N-1}wT_{1}^{-1}\cdots
T_{i-2}^{-1}\\
&  =t^{i-1-N}T_{i-1}T_{i}T_{i+1}\cdots T_{N-1}T_{i-1}wT_{1}^{-1}\cdots
T_{i-2}^{-1}\\
&  =t^{i-1-N}T_{i-1}T_{i}T_{i+1}\cdots T_{N-1}wT_{i}T_{1}^{-1}\cdots
T_{i-2}^{-1}=\xi_{i-1}T_{i}.
\end{align*}
The analogous argument as in the previous lemma shows $T_{i}\xi_{j}=\xi
_{j}T_{i}$ for $j<i-1.$Suppose $j>i+1$ then%
\begin{align*}
T_{i}\xi_{j}  &  =t^{j-N}T_{i}T_{j}T_{j+1}\cdots T_{N-1}wT_{1}^{-1}\cdots
T_{j-1}^{-1}=t^{j-N}T_{j}T_{j+1}\cdots T_{N-1}T_{i}wT_{1}^{-1}\cdots
T_{j-1}^{-1}\\
&  =t^{j-N}T_{j}T_{j+1}\cdots T_{N-1}wT_{i+1}T_{1}^{-1}\cdots T_{j-1}^{-1}\\
&  =t^{j-N}T_{j}\cdots T_{N-1}wT_{1}^{-1}\cdots T_{i-1}T_{i-2}^{-1}%
T_{i-1}^{-1}\cdots T_{j-1}^{-1}.
\end{align*}
The modified braid relations $aba=bab\Leftrightarrow ab^{-1}a^{-1}%
=b^{-1}a^{-1}b$ imply $T_{i+1}T_{i}^{-1}T_{i+1}^{-1}=T_{i}^{-1}T_{i+1}%
^{-1}T_{i`}$ and thus $T_{i}\xi_{j}=\xi_{j}T_{i}$. As before
\begin{align*}
T_{i}\xi_{i}  &  =t^{-1}T_{i}^{2}\xi_{i+1}T_{i}=t^{-1}\left\{  \left(
t-1\right)  T_{i}+t\right\}  \xi_{i+1}T_{i}=\left(  t-1\right)  \xi_{i}%
+\xi_{i+1}T_{i}.\\
\xi_{i}T_{i}  &  =\left(  t-1\right)  \xi_{i}+T_{i}\xi_{i+1}.
\end{align*}

\end{proof}

Polynomials are spanned by monomials $x^{\alpha}=\prod\limits_{i=1}^{N}%
x_{i}^{\alpha_{i}},\alpha\in\mathbb{N}_{0}^{N}$. For $\alpha\in\mathbb{N}%
_{0}^{N}$ set $s_{i}\alpha=\left(  \alpha_{1},\ldots,\overset{i}{\alpha_{i+1}%
},\overset{i+1}{\alpha_{i}},\ldots\right)  $ for $1\leq i<N$, and $\left\vert
\alpha\right\vert =\sum_{j=1}^{N}\alpha_{j}$ (the degree of $x^{\alpha}$). Let
$\mathbb{N}_{0}^{N,+}=\left\{  \alpha\in\mathbb{N}_{0}^{N}:\alpha_{1}%
\geq\alpha_{2}\geq\ldots\geq\alpha_{N}\right\}  $, the set of partitions of
length $\leq N$. Let $\alpha^{+}$ denote the nonincreasing rearrangement of
$\alpha$ (thus $\alpha^{+}\in\mathbb{N}_{0}^{N,+}$). There is a partial order
on $\mathbb{N}_{0}^{N}$%
\begin{align*}
\alpha &  \prec\beta\Longleftrightarrow\sum_{j=1}^{i}\alpha_{j}\leq\sum
_{j=1}^{i}\beta_{j},~1\leq i\leq N,~\alpha\neq\beta\text{,}\\
\alpha\lhd\beta &  \Longleftrightarrow\left(  \left\vert \alpha\right\vert
=\left\vert \beta\right\vert \right)  \wedge\left[  \left(  \alpha^{+}%
\prec\beta^{+}\right)  \vee\left(  \alpha^{+}=\beta^{+}\wedge\alpha\prec
\beta\right)  \right]  \text{,}%
\end{align*}
and a rank function ($1\leq i\leq N$)%
\[
r_{\alpha}\left(  i\right)  :=\#\left\{  j:\alpha_{j}>\alpha_{i}\right\}
+\#\left\{  j:1\leq j\leq i,\alpha_{j}=\alpha_{i}\right\}  \text{.}%
\]
Note $\alpha_{i}=\alpha_{r_{\alpha}\left(  i\right)  }^{+}$.

\subsection{Nonsymmetric Macdonald Polynomials}

The key fact about the Cherednik operators is the triangular property (see
\cite{BF1997})%
\begin{equation}
\xi_{i}x^{\alpha}=q^{\alpha_{i}}t^{N-r_{\alpha}\left(  i\right)  }x^{\alpha
}+\sum_{\beta\vartriangleleft\alpha}c_{\alpha,\beta}\left(  q,t\right)
x^{\beta}, \label{triangXi}%
\end{equation}
where the coefficients $c_{\alpha,\beta}\left(  q,t\right)  $ are polynomials
in $q,t$. For generic $\left(  q,t\right)  $ (this means $q^{m}t^{n}\neq1,0$
for $m\geq0$ and $1\leq n\leq N$) there is a basis of $\mathcal{P}$, for
$\alpha\in\mathbb{N}_{0}^{N}$%
\[
M_{\alpha}\left(  x\right)  =q^{b\left(  \alpha\right)  }t^{e\left(
\alpha\right)  }x^{\alpha}+\sum_{\beta\vartriangleleft\alpha}A_{\alpha,\beta
}\left(  q,t\right)  x^{\beta}%
\]
(where $A_{\alpha,\beta}\left(  q,t\right)  $ is a rational function of
$\left(  q,t\right)  $ with no poles when $\left(  q,t\right)  $ is generic)
and for $1\leq i\leq N$%
\[
\xi_{i}M_{\alpha}=q^{\alpha_{i}}t^{N-r_{\alpha}\left(  i\right)  }M_{\alpha}.
\]
The exponents are $b\left(  \alpha\right)  =\frac{1}{2}\sum_{i=1}^{N}%
\alpha_{i}\left(  \alpha_{i}-1\right)  $ and $e\left(  \alpha\right)
=\sum_{i=1}^{N}\alpha_{i}^{+}\left(  N-2i+1\right)  -\mathrm{inv}\left(
\alpha\right)  ,$ with $\mathrm{inv}\left(  \alpha\right)  :=\#\left\{
\left(  i,j\right)  :1\leq i<j\leq N,\alpha_{i}<\alpha_{j}\right\}  $ ; there
is an equivalent formula:%
\[
e\left(  \alpha\right)  =\frac{1}{2}\sum_{1\leq i<j\leq N}\left(  \left\vert
\alpha_{i}-\alpha_{j}\right\vert +\left\vert \alpha_{i}-\alpha_{j}%
+1\right\vert -1\right)  .
\]
These powers arise from the Yang-Baxter graph method of constructing the
$M_{\alpha}$, and are not actually needed here. The \textit{spectral vector}
of $M_{\alpha}$ is $\left[  \zeta_{\alpha}\left(  i\right)  \right]
_{i=1}^{N}$ with $\zeta_{\alpha}\left(  i\right)  =q^{\alpha_{i}%
}t^{N-r_{\alpha}\left(  i\right)  }$. We will need the formulas for the action
of $T_{i}$ on $M_{\alpha}$. Suppose $\alpha_{i}<\alpha_{i+1}$ and
$z=\zeta_{\alpha}\left(  i+1\right)  /\zeta_{\alpha}\left(  i\right)
=q^{\alpha_{i+1}-\alpha_{i}}t^{r_{\alpha}\left(  i\right)  -r_{\alpha}\left(
i+1\right)  }$ then%
\begin{align}
T_{i}M_{\alpha}  &  =M_{s_{i}\alpha}-\frac{1-t}{1-z}M_{\alpha},\label{TMi}\\
T_{i}M_{s_{i}\alpha}  &  =\frac{\left(  1-zt\right)  \left(  t-z\right)
}{\left(  1-z\right)  ^{2}}M_{\alpha}+\frac{z\left(  1-t\right)  }{\left(
1-z\right)  }M_{s_{i}\alpha}. \label{TMsi}%
\end{align}
If $\alpha_{i}=\alpha_{i+1}$ then $T_{i}M_{\alpha}=tM_{\alpha}$. The quadratic
relation appears as%
\[
\left(  T_{i}+\frac{1-t}{1-z}\right)  \left(  T_{i}-\frac{z\left(  1-t\right)
}{1-z}\right)  =\frac{\left(  1-zt\right)  \left(  t-z\right)  }{\left(
1-z\right)  ^{2}}.
\]

\subsection{Action of $T_{i}$ on polynomials and $\vartriangleright$-maximal
terms}

The following are routine computations:

\begin{lemma}
\label{Txgm}Suppose $\gamma\in\mathbb{N}_{0}^{N}$ and $1\leq i<N$ . Set
$x^{\prime}=\prod_{j\neq i,i+1}x^{\gamma_{j}}$. Then\newline(1) $\gamma
_{i}>\gamma_{i+1}+1$ implies $T_{i}x^{\gamma}=\left(  1-t\right)  x^{\prime
}\sum\limits_{j=0}^{\gamma_{i}-\gamma_{i+1}-1}x_{i}^{\gamma_{i}-j-1}%
x_{i+1}^{\gamma_{i+1}+j+1}+tx^{s_{i}\gamma}$;\newline(2) $\gamma_{i}%
=\gamma_{i+1}+1$ implies $T_{i}x^{\gamma}=x^{s_{i}\gamma};$\newline(3)
$\gamma_{i}=\gamma_{i+1}$ implies $T_{i}x^{\gamma}=tx^{\gamma}$;\newline(4)
$\gamma_{i}=\gamma_{i+1}-1$ implies $T_{i}x^{\gamma}=tx^{s_{i}\gamma}+\left(
t-1\right)  x^{\gamma};$\newline(5) $\gamma_{i}<\gamma_{i+1}-1$ implies
$T_{i}x^{\gamma}=\left(  t-1\right)  x^{\prime}\sum\limits_{j=0}^{\gamma
_{i+1}-\gamma_{i}-1}x_{i}^{\gamma_{i}+j}x_{i+1}^{\gamma_{i+1}-j}%
+tx^{s_{i}\gamma}$.
\end{lemma}

\begin{lemma}
Suppose $\lambda\in\mathbb{N}_{0}^{N,+},\lambda_{i}>\lambda_{j}+1$ ($i>j$) and
$1\leq s<\lambda_{i}-\lambda_{j}$, $\mu\in\mathbb{N}_{0}^{N}$ such that
$\mu_{k}=\lambda_{k}$ for $k\neq i,j$, $\mu_{i}=\lambda_{i}-s,\mu_{j}%
=\lambda_{j}+s$ then $\lambda\succ\mu^{+}$.
\end{lemma}

(The proof is left as an exercise.)

In (1) above let $\alpha_{k}=\gamma_{k}$ for $k\neq i,i+1$ and $\alpha
_{i}=\gamma_{i}-j-1,\alpha_{i+1}=\gamma_{i+1}+j+1$ with $1\leq j+1<\gamma
_{i}-\gamma_{i+1}$ then the Lemma with $\lambda=\gamma^{+}$ and $\mu
^{+}=\alpha^{+}$ shows $\gamma^{+}\succ a^{+}$ (the other term in (1) is
$x^{s_{i}\gamma}$ and $\gamma\succ s_{i}\gamma$). Similarly in (5) let
$\alpha_{i}=\gamma_{i}+j,\alpha_{i+1}=\gamma_{i+1}-j$ with $1\leq j\leq
\gamma_{i+1}-\gamma_{i}-1$, thus $\gamma^{+}\succ a^{+}$ (the other term in
(5) for $j=0$ is $x^{\gamma}$ and $s_{i}\gamma\succ\gamma$.

\begin{proposition}
\label{maxalf}Suppose $\alpha$ is $\vartriangleright$-maximal in
$p=\sum_{\delta}c_{\delta}x^{\delta}$ (a homogeneous polynomial, $\left\vert
\delta\right\vert =\left\vert \alpha\right\vert $), that is, $c_{\alpha}\neq0$
and if some $\delta\trianglerighteq\alpha$ with $c_{\delta}\neq0$ then
$\delta=\alpha$. Furthermore suppose $\alpha_{i+1}>\alpha_{i}$ for some i and
$x^{\beta}$ with $\beta\vartriangleright s_{i}\alpha$ appears in $\left(
T_{i}+c\right)  p$ then $\beta^{+}=\alpha^{+}$ and $\beta\succ s_{i}\alpha$.
\end{proposition}

\begin{proof}
Suppose $x^{\beta}$ appears in $T_{i}x^{\gamma}$ (with $c_{\gamma}\neq0$) in
one of the five cases of Lemma \ref{Txgm} and $\beta^{+}\succ\left(
s_{i}\alpha\right)  ^{+}=\alpha^{+}.$ Every term satisfies $\gamma^{+}%
\succ\beta^{+}$ or $\gamma^{+}=\beta^{+}$ but then $\gamma^{+}\succeq\beta
^{+}\succ\alpha^{+}$ and $\gamma\vartriangleright\alpha$, a contradiction.
Suppose $\beta^{+}=\alpha^{+}$ then $\beta\vartriangleright s_{i}\alpha$
implies $\beta\succ s_{i}\alpha$.
\end{proof}

\begin{corollary}
If $\alpha$ is $\vartriangleright$-maximal in $p=\sum_{\delta}c_{\delta
}x^{\delta}$ and $x^{\beta}$ appears in $\left(  T_{i}+c\right)  p$ with
$\beta\trianglerighteq s_{i}\alpha$ then either $\beta=s_{i}\alpha$ or
$\beta^{+}=\alpha^{+}$ and $\beta\succ s_{i}\alpha$ with $\beta=s_{i}\gamma$
where $x^{\gamma}$ appears in $p$.
\end{corollary}

\begin{proof}
If $\beta$ occurs in case (1) or case (5) of Lemma \ref{Txgm} and $\beta
\neq\gamma,s_{i}\gamma$ (for $x^{\gamma}$ appearing in $p$) then $\gamma
^{+}\succ\beta^{+}\succ s_{i}\alpha\succ\alpha$ which violates the
$\vartriangleright$-maximality of $\alpha$ this leaves only $\beta=s_{i}%
\gamma$.
\end{proof}

Note $\beta=s_{i}\gamma$ does not imply $s_{i}\beta\succ\alpha$, for example
let $\beta=\left(  4,1,3,2\right)  $ and $s_{1}\alpha=\left(  3,2,1,4\right)
$ then $\beta\succ s_{1}\alpha$ but $s_{1}\beta=\left(  1,4,3,2\right)  $ and
$\alpha=\left(  2,3,1,4\right)  $ are not $\vartriangleright$-comparable.

\subsection{Irreducible representations of the Hecke algebra}

Irreducible representations of $\mathcal{H}_{N}\left(  t\right)  $ are indexed
by partitions of $N$ (for background see Dipper and James \cite{DJ1986}).
Given a partition $\tau\in\mathbb{N}_{0}^{N,+}$ with $\left\vert
\tau\right\vert =N$ there is a \textit{Ferrers diagram}: boxes at $\left(
i,j\right)  $ with $1\leq i\leq\ell\left(  \tau\right)  =\max\left\{
j:\tau_{j}>0\right\}  $ and $1\leq j\leq\tau_{i}$. The module is spanned by
reverse standard Young tableaux (abbr. RSYT) of shape $\tau$ (denoted
$\mathcal{Y}_{\tau})$: the numbers $1,\ldots,N$ are inserted into the Ferrers
diagram so that the entries in each row and in each column are decreasing. The
module $\mathrm{span}_{\mathbb{K}}\left\{  Y:Y\in\mathcal{Y}_{\tau}\right\}  $
is said to be of isotype $\tau$. If $k$ is in cell $\left(  i,j\right)  $ of
RSYT $Y$ (denoted $Y\left[  i,j\right]  =k$) then the \textit{content}
$c\left(  k,Y\right)  :=j-i$; the \textit{content vector} $\left[  c\left(
k,Y\right)  \right]  _{k=1}^{N}$ determines $Y$ uniquely. The action of
$\mathcal{H}_{N}\left(  t\right)  $ is specified by the formulas for $T_{i}Y$:

\begin{itemize}
\item if $c\left(  i,Y\right)  -c\left(  i+1,Y\right)  =1$ then $T_{i}Y=tY$;

\item if $c\left(  i,Y\right)  -c\left(  i+1,Y\right)  =-1$ then $T_{i}Y=-Y;$

\item if $\left\vert c\left(  i,Y\right)  -c\left(  i+1,Y\right)  \right\vert
\geq2$ then let $Y^{\left(  i\right)  }$ denote the RSYT obtained by
interchanging $i$ and $i+1$ in $Y$ and set $z=t^{c\left(  i+1,Y\right)
-c\left(  i,Y\right)  }$: if $c\left(  i,Y\right)  -c\left(  i+1,Y\right)
\geq2$ then%
\[
T_{i}Y=Y^{\left(  i\right)  }-\frac{1-t}{1-z}Y;
\]
if $c\left(  i,Y\right)  -c\left(  i+1,Y\right)  \leq-2$ then%
\[
T_{i}Y=\frac{\left(  1-zt\right)  \left(  t-z\right)  }{\left(  1-z\right)
^{2}}Y^{\left(  i\right)  }-\frac{1-t}{1-z}Y.
\]

\end{itemize}

From these relations it follows that $\omega_{i}Y=t^{c\left(  i,Y\right)  }Y$
for $1\leq i\leq N$. Call the vector $\left[  t^{c\left(  i,Y\right)
}\right]  _{i=1}^{N}$ the $t$-exponential content vector of $Y$, or the
$t^{C}$-vector for short. Note $c\left(  N,Y\right)  =0$ always and
$\omega_{N}:=1.$

So if one finds a simultaneous eigenfunction of $\left\{  \omega_{i}\right\}
$ then the eigenvalues determine an RSYT and the isotype (partition) of an
irreducible. representation.

\subsection{Singular parameters}

For integers $m$ and $n$ such that $m\geq1$ and $2\leq n\leq N$ we consider
\textit{singular} parameters $\left(  q,t\right)  $ satisfying $q^{m}t^{n}=1$
with the property that if $q^{a}t^{b}=1$ then $a=rm,b=rn$ for some
$r\in\mathbb{Z}$.

\begin{definition}
\label{qmtndef}Let $g=\gcd(m,n)$ and let ${z=\exp\left(  \frac{2\pi
\mathrm{i}k}{m}\right)  }$ with $\gcd\left(  k,g\right)  =1$, that is,
$z^{m/g}$ is a primitive $g^{\text{th}}$ root of unity. If $g=1$ then set
$z=1.$ Define $\varpi:=\left(  q,t\right)  =\left(  zu^{-n/g},u^{m/g}\right)
$ where $u$ is not a root of unity and $u\neq0$.
\end{definition}

\begin{lemma}
\label{qmtn}If $q^{a}t^{b}|_{\varpi}=1$ for some integers $a,b$ then
$a=rm,b=rn$ for some $r\in\mathbb{Z}$.
\end{lemma}

\begin{proof}
By hypothesis $z^{a}u^{-an/g+bm/g}=1$ and, since $u$ is not a root of unity,
${-a\frac{n}{g}+b\frac{m}{g}=0}$. From ${\gcd\left(  \frac{n}{g},\frac{m}%
{g}\right)  =1}$, it follows that $a=p^{\prime}\frac{m}{g}$ and $b=p^{\prime
}\frac{n}{g}$, for some $p^{\prime}\in\mathbb{Z}$. Thus ${1=z^{a}=\exp\left(
\frac{2\pi\mathrm{i}k}{m}\frac{mp^{\prime}}{g}\right)  =\exp\left(  \frac
{2\pi\mathrm{i}k}{g}p^{\prime}\right)  }$. Moreover, since $\gcd\left(
k,g\right)  =1$, $p^{\prime}=pg$ with $p\in\mathbb{Z}$. Hence $a=pm$ and
$b=pn$.
\end{proof}

In fact, to describe all the possibilities for $\varpi$, it suffices to let
$1\leq k<g$. To be precise, $\varpi$ is not a single point but a variety in
$\left(  \mathbb{C}\backslash\left\{  0\right\}  \right)  ^{2}$.

\section{\label{necessary}Necessary Conditions for Singular Polynomials}

By using the degree-lowering ($q$-Dunkl) operators defined by Baker and
Forrester \cite{BF1997} we find another characterization of singular polynomials.

\begin{definition}
Suppose $p\in\mathcal{P}$ then%
\begin{align*}
D_{N}p\left(  x\right)   &  :=\frac{1}{x_{N}}\left(  1-\xi_{N}\right)
p\left(  x\right)  ,\\
D_{i}p\left(  x\right)   &  :=\dfrac{1}{t}T_{i}D_{i+1}T_{i}p\left(  x\right)
,~i<N.
\end{align*}

\end{definition}

\begin{proposition}
A polynomial $p$ is singular if and only if $D_{i}p=0$ for $1\leq i\leq N$.
\end{proposition}

\begin{proof}
The proof is by downward induction on $i$. Since $\omega_{N}=1$ it follows
that $D_{N}p=0$ iff $\xi_{N}p=p=\omega_{N}p$. Suppose that $D_{i}p=0$ iff
$\xi_{i}p=\omega_{i}p$ for all $p$ and $k\leq i\leq N.$ Then $D_{k-1}p=0$ iff
$t^{-1}T_{k-1}D_{k}T_{k-1}p=0$ iff $D_{k}T_{k-1}p=0$ iff $\xi_{k}%
T_{k-1}p=\omega_{k}T_{k-1}p$ iff $t^{-1}T_{k-1}\xi_{k}T_{k-1}p=t^{-1}%
T_{k-1}\omega_{k}T_{k-1}p$.
\end{proof}

First we show that any singular polynomial generates an $\mathcal{H}%
_{N}\left(  t\right)  $-module consisting of singular polynomials. This allows
the use of the representation theory of $\mathcal{H}_{N}\left(  t\right)  $.

\begin{proposition}
Suppose $p$ is singular and $1\leq i<N$, then $T_{i}p$ is singular.
\end{proposition}

\begin{proof}
The commutation relations from Lemmas \ref{comrel1} and \ref{comrel2} are
used. Suppose $j<i$ or $j>i+1$ then $\xi_{j}T_{i}p=T_{i}\xi_{j}p=T_{i}%
\omega_{j}p=\omega_{j}T_{i}p$. Case $j=i$:%
\begin{align*}
\xi_{i}T_{i}p  &  =\left\{  \left(  t-1\right)  \xi_{i}+T_{i}\xi
_{i+1}\right\}  p=\left(  t-1\right)  \omega_{i}p+T_{i}\omega_{i+1}p\\
&  =\left\{  \left(  t-1\right)  \omega_{i}+T_{i}\omega_{i+1}\right\}
p=\omega_{i}T_{i}p.
\end{align*}
Case $j=i+1$%
\begin{align*}
\xi_{i+1}T_{i}p  &  =\left\{  T_{i}\xi_{i}-\left(  t-1\right)  \xi
_{i}\right\}  p=T_{i}\omega_{i}p-\left(  t-1\right)  \omega_{i}p\\
&  =\left\{  T_{i}\omega_{i}-\left(  t-1\right)  \omega_{i}\right\}
p=\omega_{i+1}T_{i}p.
\end{align*}

\end{proof}

\begin{proposition}
Suppose $p$ is singular then $\mathcal{M}=\mathcal{H}_{N}\left(  t\right)  p$
is a linear space of singular polynomials, and it is closed under the actions
of $\xi_{i},\omega_{i}.$for $1\leq i\leq N$, and $w.$
\end{proposition}

\begin{proof}
By definition of $\omega_{i}$ we see that $f\in\mathcal{M}$ implies
$\omega_{i}f\in\mathcal{M}$, and by definition $\xi_{i}f=\omega_{i}%
f\in\mathcal{M}$. Also
\begin{align*}
\xi_{1}p  &  =T_{1}T_{2}\cdots T_{N-1}wp\\
&  =\omega_{1}p=t^{1-N}T_{1}T_{2}\ldots\cdots T_{N-1}T_{N-1}T_{N-2}\cdots
T_{1}p
\end{align*}
thus $wp=t^{1-N}T_{N-1}T_{N-2}\cdots T_{1}p$.
\end{proof}

Note that $\mathcal{M}$ is also a module of the affine Hecke algebra. By the
representation theory of $\mathcal{H}_{N}\left(  t\right)  $ the module has a
basis of $\left\{  \omega_{i}\right\}  $-simultaneous eigenfunctions and by
definition these are $\left\{  \xi_{i}\right\}  $-simultaneous eigenfunctions
- note we are not claiming they are specializations of nonsymmetric Macdonald
polynomials at $\varpi$. Suppose $f$ is such an eigenfunction and let $\alpha$
be $\vartriangleright$-maximal in the expression $f\left(  x\right)
=\sum_{\beta}c_{\beta}x^{\beta}$. Then $\xi_{i}f=q^{\alpha_{i}}t^{N-r_{\alpha
}\left(  i\right)  }f$ because by the triangularity property of $\xi_{i}$ (see
(\ref{triangXi})) $x^{\alpha}$ can only appear in $\xi_{i}f$ in the term
$\xi_{i}x^{\alpha}$. Furthermore $\xi_{i}f=\omega_{i}f$ implies $q^{\alpha
_{i}}t^{N-r_{\alpha}\left(  i\right)  }=t^{c\left(  i,Y\right)  }$ for some
RSYT $Y$, at $\varpi$. As well we can conclude $\alpha_{i}=mr,N-r_{\alpha
}\left(  i\right)  -c\left(  i,Y\right)  =nr$ for some $r\in\mathbb{N}$ (Lemma
\ref{qmtn}). The next step is to produce a simultaneous eigenfunction which
has a $\vartriangleright$-maximal term $x^{\lambda}$ with $\lambda
\in\mathbb{N}_{0}^{N,+}$.

\begin{proposition}
There exists $f\in\mathcal{M}$ which is a simultaneous $\left\{  \omega
_{i}\right\}  $-eigenfunction and $f=c_{\lambda}x^{\lambda}+\sum
_{\beta\vartriangleleft\lambda}c_{\beta}x^{\beta}+\sum_{\gamma}c_{\gamma
}x^{\gamma}$ where $\gamma$ is not $\vartriangleright$-comparable to $\lambda
$, and $\lambda\in\mathbb{N}_{0}^{N,+}$.
\end{proposition}

\begin{proof}
Suppose $f=\sum c_{\alpha}x^{\alpha}$ is an eigenfunction and there is a
$\vartriangleright$-maximal $\alpha$ with $x^{\alpha}$ (i.e. $c_{\alpha}\neq
0$) appearing in $f$, and $\alpha_{i}<\alpha_{i+1}$ then $T_{i}f\neq f$ and
the coefficient of $x^{s_{i}\alpha}$ is $tc_{\alpha}$; let $\omega_{j}%
f=\mu_{j}f$ for $1\leq j\leq N$ and $\mu_{i+1}\neq\mu_{i}$ (because $c\left(
i,Y\right)  \neq c\left(  i+1,Y\right)  $ for any RSYT) so that%
\[
g:=T_{i}f+\dfrac{t-1}{\mu_{i+1}/\mu_{i}-1}f
\]
is a simultaneous eigenfunction with $\vartriangleright$-maximal $\beta$ such
that $\beta^{+}=\alpha^{+}$ and $\beta\succeq s_{i}\alpha$, (by Proposition
\ref{maxalf}) and eigenvalues \ldots$\mu_{i+1},\mu_{i}\ldots$In general this
formula could produce a zero function $g$ but this does not happen here
because the coefficient of $x^{s_{i}\alpha}$ in $g$ is not zero. Repeating
these steps eventually produces a $\vartriangleright$-maximal term
$x^{\lambda}$ with $\lambda\in\mathbb{N}_{0}^{N,+}$ (at most $\mathrm{inv}%
\left(  \alpha\right)  $ steps).
\end{proof}

At this point we have shown if there is a singular polynomial then there is a
partition $\lambda\in\mathbb{N}_{0}^{N,+}$ and an RSYT\ $Y$ such that
$q^{\lambda_{i}}t^{N-i}=t^{c\left(  i,Y\right)  }$ at $\varpi$, for $1\leq
i\leq N$. Next we determine necessary conditions on $\lambda$ for the
existence of $Y$, in other words, when $\left[  q^{\lambda_{i}}t^{N-i}\right]
_{i=1}^{N}$ at $\varpi$ is a valid $t^{C}$-vector. The equations $\lambda
_{i}=mr_{i},N-i-c\left(  i,Y\right)  =nr_{i}$ for $1\leq i\leq N$ show that
$\lambda$ can be replaced by $\frac{1}{m}\lambda$ and $\varpi$ by $qt^{n}=1$
(simply $q=t^{-n}$), also $n\lambda_{i}=N-i-c\left(  i,Y\right)  $.

The following is a restatement of the development in \cite{CD2020} with
significant differences in notation. First there is an informal discussion of
the beginning of the process of building $Y$ by placing $N,N-1,N-2,\ldots$in
possible locations and determining $\lambda_{N},\lambda_{N-1},\lambda
_{N-2},\ldots$accordingly. Abbreviate $c_{i}=c\left(  i,Y\right)  $.

Suppose $\lambda_{N-k}$ is the last nonzero entry of $\lambda$ ($\lambda
_{i}=0$ for $i>N-k$) then $k-c_{N-k}=n\lambda_{N-k}$; the entry $N-k$ in $Y$
is at $\left[  1,k+1\right]  $ or $\left[  2,1\right]  $ thus $c_{N-k}%
=k,\lambda_{N-k}=0$ (contra) or $c_{N-k}=-1,n\lambda_{N-k}=k+1$. Set
$\lambda_{N-k}=d_{1}$ and $k=nd_{1}-1.$The entry $N-k-1$ in $Y$ is in one of
$\left[  3,1\right]  ,\left[  2,2\right]  ,\left[  1,k+1\right]  $ with
contents $-2,0,k$ respectively, yielding the equations $n\lambda
_{N-k-1}=k+1-c_{N-k-1}=k-1,k+1,1=nd_{1}-2,nd_{1},1$, respectively. If $n>2$
then only $\left[  2,2\right]  $ is possible and $\lambda_{N-k-1}=d_{1}$. If
$n=2$ then $\left[  3,1\right]  ,\lambda_{N-k-1}=d_{1}+1$ and $\left[
2,2\right]  ,\lambda_{N-k-1}=d_{1}$ are possible.

\begin{theorem}
\label{singLB}There are numbers $d_{1}\geq d_{2}\geq\ldots\geq d_{L}\geq1$
such that with $\gamma_{s}:=\sum_{i=1}^{s-1}d_{i}$ and $0\leq r_{L+1}%
<N-n\gamma_{L+1}+L\leq nd_{L}-1$ the entries in row $s$ of $Y$ are
$R_{s}:=\left\{  i:n\gamma_{s}-s+1\leq N-i\leq n\gamma_{s+1}-s-1\right\}  $
for $1\leq s\leq L$,\newline$R_{L+1}=\left\{  i:n\gamma_{L+1}-L\leq N-i\leq
N-1\right\}  $ and $\lambda_{i}=\gamma_{s}$ for $i\in R_{s}$. The isotype of
$Y$ is $\tau:=\left(  nd_{1}-1,nd_{2}-1,\ldots,nd_{L}-1,r_{L+1}\right)  $.
\end{theorem}

\begin{proof}
By way of induction suppose there are numbers $d_{1}\geq d_{2}\geq\ldots\geq
d_{k-1}>0$ such that the entries in row $s$ of $Y$ are\newline\ $R_{s}%
=\left\{  i:n\gamma_{s}-s+1\leq N-i\leq n\gamma_{s+1}-s-1\right\}  $ and
$\lambda_{i}=\gamma_{s}$ for $i\in R_{s}$. Assume this has been proven for
$1\leq s<k$ and for row $k$ up to $n\gamma_{k}-k+1\leq N-i\leq n\gamma
_{k}-k+\ell$ with $\ell\leq nd_{k-1}-1$ (the length $\#R_{k-1}$ of row $k-1$).
Consider the possible locations for the next entry $p=N-\left(  n\gamma
_{k}-k+\ell+1\right)  $. The possible boxes are (1) $\left[  s,nd_{s}\right]
$ ($s<k$ and $d_{s}<d_{s-1}$ or $s=1$), (2) $\left[  k,\ell+1\right]  $, (3)
$\left[  k+1,1\right]  $ with contents $nd_{s}-s,\ell+1-k,-k$ respectively.
The equations%
\begin{align*}
n\lambda_{p} &  =N-p-c_{p}=n\gamma_{k}-k+\ell+1-c_{p}\\
n\left(  \lambda_{p}-\gamma_{k}\right)   &  =-k+\ell+1-c_{p}%
\end{align*}
must hold;\newline case (1): (note $\ell+1\leq nd_{k-1}$)%
\begin{align*}
n\left(  \lambda_{p}-\gamma_{k}\right)   &  =-k+\ell+1-nd_{s}+s\\
n\left(  \lambda_{p}-\gamma_{k}+d_{s}\right)   &  =-k+s+1+\ell\leq
-k+s+nd_{k-1}\\
n\left(  \lambda_{p}-\gamma_{k}+d_{s}-d_{k-1}\right)   &  \leq s-k<0
\end{align*}
$\lambda_{p}\geq\gamma_{k}=\lambda_{p+1}$ and $d_{s}\geq d_{k-1}$ by inductive
hypothesis, so the left side $\geq0$ and there is a contradiction.\newline
case (2):%
\begin{align*}
n\left(  \lambda_{p}-\gamma_{k}\right)   &  =-k+\ell+1-\left(  \ell
+1-k\right)  =0\\
\lambda_{p} &  =\gamma_{k}%
\end{align*}
and the inductive hypothesis is proved for $n\gamma_{k}-k+1\leq N-i\leq
n\gamma_{k}-k+\ell+1$, entries in row $k$.\newline case (3)
\[
n\left(  \lambda_{p}-\gamma_{k}\right)  =-k+\ell+1+k=\ell+1
\]
set $\ell=nd_{k}-1$ and $\gamma_{k+1}=\gamma_{k}+d_{k},\lambda_{p}%
=\gamma_{k+1}$. The inductive step has been proven for $k$ and for $k+1$ with
$Y\left[  k+1,1\right]  =N-n\gamma_{k+1}+k$. By induction this uses up all the
entries. Let row $L+1$ be the last row of $Y$ and of length $r_{L+1}$, then
$N=\sum_{i=1}^{L}\left(  nd_{i}-1\right)  +r_{L+1}$ and $r_{L+1}\leq nd_{L}-1$.
\end{proof}

\begin{corollary}
\label{omxiEv}Suppose $\varpi=\left(  q,t\right)  $ as in Definition
\ref{qmtndef} and $p$ is singular. Then $\mathcal{H}_{N}\left(  t\right)  p$
contains a $\left\{  \omega_{i},\xi_{i}\right\}  $ simultaneous eigenfunction
$f=c_{\lambda}x^{\lambda}+\sum_{\beta\vartriangleleft\lambda}c_{\beta}%
x^{\beta}+\sum_{\gamma}c_{\gamma}x^{\gamma}$ with $\gamma$ not
$\vartriangleright$-comparable to $\lambda$ so that $\lambda_{i}=m\gamma_{s}$
if $i\in R_{s}$, in the notation of the Theorem.
\end{corollary}

We have shown if $\alpha$ is $\vartriangleright$-maximal in a simultaneous
$\left\{  \omega_{i},\xi_{i}\right\}  $ eigenfunction then there is an
eigenfunction in which $\alpha^{+}$ is $\vartriangleright$-maximal. Now the
eigenvalues are determined by $Y$ and it follows that $\alpha^{+}=\lambda$ as
constructed above. Hence each term $x^{\gamma}$ in an eigenfunction satisfies
$\gamma\trianglelefteq\lambda$. (Suppose at some stage $\gamma$ is
$\vartriangleright$-maximal then there is a simultaneous eigenfunction with
$\gamma^{+}$ being $\vartriangleright$-maximal and the construction produces
an RSYT of the same isotype $\tau$ and the numbers $N,N-1,\ldots$ are entered
row-by-row forcing $\gamma^{+}=\lambda$.)

\begin{theorem}
(\cite{CD2020}) In the notation of Theorem \ref{singLB} if $d_{i}=1$ for
$i\geq2$ then $M_{\lambda}\left(  x\right)  $ specialized to $\varpi$ has no
poles and is singular. The module $\mathcal{H}_{N}\left(  t\right)
M_{\lambda}$ is spanned by $M_{\alpha\left(  Y\right)  }$ where $Y\in
\mathcal{Y}_{\tau}$,$\tau=\left(  nd_{1}-1,\left(  n-1\right)  ^{L-1}%
,r_{L+1}\right)  $ and $\alpha\left(  Y\right)  _{i}=m\left(  d_{1}%
+s-2\right)  $ if $Y\left[  s,k\right]  =i$ for $s\geq2$ and some $k$,
otherwise $(Y\left[  1,k\right]  =i$) $\alpha\left(  Y\right)  _{i}=0$.
\end{theorem}

The Ferrers diagram of $\lambda$ (from Theorem \ref{singLB}) is called a
quasi-staircase, the shape suggested when French notation with row 1 on the
bottom is used.

We have reached the main purpose of this paper: to show there are no other
singular polynomials.

\section{\label{restrict}Restrictions}

In this section we show that the desired nonexistence result can be reduced to
the simpler two-row situation.

Suppose $\alpha\in\mathbb{N}_{0}^{N}$ and $r_{\alpha}\left(  1\right)  =1$
(that is, $\alpha_{i}\leq\alpha_{1}$ for all $i$). Let $\alpha^{\prime
}=\left(  \alpha_{2},\ldots,\alpha_{N}\right)  $ and $Y^{\prime}%
=Y\backslash\left\{  1\right\}  $ (the RSYT where the entry $1$ is deleted)
and $f$ satisfies $\xi_{i}f=q^{\alpha_{i}}t^{N-r_{\alpha}\left(  i\right)  }%
f$, at $\varpi$. First we will show that $f_{\alpha^{\prime}}:=\mathrm{coeff}%
\left(  x_{1}^{\alpha_{1}},f\right)  $ is an eigenfunction of $\xi
_{i}^{^{\prime}}$ with eigenvalue $q^{\alpha_{i}}t^{N-r_{\alpha}\left(
i\right)  }$ for $2\leq i\leq N$ where
\begin{align*}
w^{\prime}p\left(  x\right)   &  :=p\left(  qx_{N},x_{2},x_{3},\ldots
,x_{N-1}\right)  ,\\
\xi_{i}^{\prime}p\left(  x\right)   &  :=t^{i-2}T_{i}T_{i+1}\cdots
T_{N-1}w^{\prime}T_{2}^{-1}\cdots T_{i-1}^{-1}p\left(  x\right)
\end{align*}

\begin{lemma}
Let $f=x_{1}^{\alpha_{1}}x_{2}^{\alpha_{2}}p\left(  x_{3},\ldots,x_{N}\right)
$ with $\alpha_{1}\geq\alpha_{2}$ then%
\[
\mathrm{coeff}\left(  x_{1}^{\alpha_{1}},wT_{1}^{-1}f\right)  =t^{-1}%
w^{\prime}\mathrm{coeff}\left(  x_{1}^{\alpha_{1}},f\right)  .
\]

\end{lemma}

\begin{proof}
By definition%
\begin{align*}
T_{1}^{-1}f  &  =\frac{1-t}{t}x_{1}\frac{f\left(  x\right)  -f\left(
xs_{1}\right)  }{x_{1}-x_{2}}+t^{-1}f\left(  xs_{1}\right) \\
&  =\frac{1-t}{t}x_{1}^{1+\alpha_{2}}x_{2}^{\alpha_{2}}\frac{x_{1}^{\alpha
_{1}-\alpha_{2}}-x_{2}^{\alpha_{1}-\alpha_{2}}}{x_{1}-x_{2}}p+t^{-1}%
x_{1}^{\alpha_{2}}x_{2}^{\alpha_{1}}p\left(  x_{3},\ldots,x_{N}\right) \\
&  =\frac{1-t}{t}\sum_{i=0}^{\alpha_{1}-\alpha_{2}-1}x_{1}^{\alpha_{1}-i}%
x_{2}^{\alpha_{2}+i}p+t^{-1}x_{1}^{\alpha_{2}}x_{2}^{\alpha_{1}}p\left(
x_{3},\ldots,x_{N}\right)
\end{align*}
then%
\begin{align*}
wT_{1}^{-1}f  &  =\frac{1-t}{t}\sum_{i=0}^{\alpha_{1}-\alpha_{2}-1}\left(
qx_{N}\right)  ^{\alpha_{1}-i}x_{1}^{\alpha_{2}+i}p\left(  x_{2},x_{3}%
,\ldots,x_{N-1}\right) \\
&  +x_{1}^{\alpha_{1}}\left(  qx_{N}\right)  ^{\alpha_{2}}t^{-1}p\left(
x_{2},x_{3},\ldots,x_{N-1}\right)  .
\end{align*}
The highest power of $x_{1}$ in the first term is $\alpha_{1}-1$ thus%
\[
\mathrm{coeff}\left(  x_{1}^{n},wT_{1}^{-1}f\right)  =\left(  qx_{N}\right)
^{\alpha_{2}}t^{-1}p\left(  x_{2},x_{3},\ldots,x_{N-1}\right)
\]
and the right hand side is $t^{-1}w^{\prime}x_{2}^{\alpha_{2}}p\left(
x_{3},\ldots,x_{N}\right)  $.
\end{proof}

Let $\pi_{n}f:=\mathrm{coeff}\left(  x_{1}^{n},f\right)  $.

\begin{theorem}
\label{xiPieq}Suppose $f=\sum_{\alpha}c_{\alpha}x^{\alpha}$ with $\max
_{i}\alpha_{i}=n$ then $\pi_{n}\xi_{i}f=\xi_{i}^{\prime}\pi_{n}f$ for $2\leq
i\leq N$.
\end{theorem}

\begin{proof}
Let $i>1$ then%
\begin{align*}
\pi_{n}\xi_{i}f  &  =t^{i-1}\pi_{n}T_{i}T_{i+1}\cdots T_{N-1}wT_{1}^{-1}%
T_{2}^{-1}\cdots T_{i-1}^{-1}f\left(  x\right) \\
&  =t^{i-1}T_{i}T_{i+1}\cdots T_{N-1}\pi_{n}wT_{1}^{-1}T_{2}^{-1}\cdots
T_{i-1}^{-1}f\left(  x\right) \\
&  =t^{i-2}T_{i}T_{i+1}\cdots T_{N-1}w^{\prime}\pi_{n}T_{2}^{-1}\cdots
T_{i-1}^{-1}f\left(  x\right) \\
&  =\ t^{i-2}T_{i}T_{i+1}\cdots T_{N-1}w^{\prime}T_{2}^{-1}\cdots T_{i-1}%
^{-1}\pi_{n}f\left(  x\right) \\
&  =\xi_{i}^{\prime}\pi_{n}f;
\end{align*}
this uses the Lemma and the fact that $\xi_{i}f$ and $T_{2}^{-1}\cdots
T_{i-1}^{-1}f$ are sums of monomials $x^{\beta}$ with $\beta_{j}\leq n$ for
$j\geq1$ (properties of the order $\vartriangleright$ and of $T_{j}^{-1}$). If
$i=2$ then the empty product $T_{2}^{-1}\cdots T_{i-1}^{-1}$ reduces to $1$.
\end{proof}

Suppose $\alpha,\beta\in\mathbb{N}_{0}^{N-1}$ (indexed $2\leq i\leq N$) and
$\left\vert \alpha\right\vert =\left\vert \beta\right\vert ,$ set
$\alpha^{\prime}=\left(  n,\alpha\right)  ,\beta^{\prime}:=\left(
n,\beta\right)  $ (so that $\left\vert \alpha^{\prime}\right\vert =\left\vert
\beta^{\prime}\right\vert ).$

\begin{lemma}
\label{resagtb}Suppose $\max_{i}\alpha_{i}\leq n$ and $\max_{i}\beta_{i}\leq
n$ then $\alpha^{\prime+}=\left(  n,\alpha^{+}\right)  ,\beta^{\prime
+}=\left(  n,\beta^{+}\right)  $ and $\alpha^{\prime}\succ\beta^{\prime}$ iff
$\alpha\succ\beta$, $\alpha^{\prime}\vartriangleright\beta^{\prime}$ iff
$\alpha\vartriangleright\beta$.
\end{lemma}

\begin{proof}
By hypothesis $\left(  \alpha^{\prime+}\right)  _{1}=n$ and $\alpha^{\prime
+}=\left(  n,\alpha^{+}\right)  $, similarly $\beta^{\prime+}=\left(
n,\beta^{+}\right)  $. Furthermore%
\begin{align*}
\alpha^{\prime}  &  \succ\beta^{\prime}\Longleftrightarrow n+\sum_{j=2}%
^{i}\alpha_{j}\geq n+\sum_{j=2}^{i}\beta_{j}~\forall i\geq2\\
&  \Longleftrightarrow\alpha\succ\beta
\end{align*}
Then%
\begin{align*}
\alpha\vartriangleright\beta &  \Longleftrightarrow\left(  \alpha^{+}%
\succ\beta^{+}\right)  \vee\left(  \alpha^{+}=\beta^{+}\wedge\alpha\succ
\beta\right) \\
\alpha^{\prime}\vartriangleright\beta^{\prime}  &  \Longleftrightarrow\left(
\alpha^{\prime+}\succ\beta^{\prime+}\right)  \vee\left(  \alpha^{\prime
+}=\beta^{\prime+}\wedge\alpha^{\prime}\succ\beta^{\prime}\right)
\end{align*}
and $\alpha\vartriangleright\beta\Longleftrightarrow\alpha^{\prime
}\vartriangleright\beta^{\prime}$.
\end{proof}

\begin{proposition}
\label{ResXiPi}Let $f$ be the $\left\{  \omega_{i},\xi_{i}\right\}  $
simultaneous eigenfunction from Corollary \ref{omxiEv} with eigenvalues
$q^{\lambda_{i}}t^{N-i}=t^{c\left(  i,Y\right)  }$ at $q^{m}t^{n}=1$ for
$1\leq i\leq N$. Then $\pi_{\lambda_{1}}f$ is a nonzero $\left\{  \omega
_{i},\xi_{i}^{\prime}:i\geq2\right\}  $ simultaneous eigenfunction with the
same eigenvalues as $f$ for $i\geq2$ with $c\left(  i,Y\right)  =c\left(
i,Y\backslash\left\{  1\right\}  \right)  .$Here $Y\backslash\left\{
1\right\}  $ is the RSYT obtained by removing the box containing $1$ from $Y$.
\end{proposition}

\begin{proof}
We showed that each term $x^{\alpha}$ appearing in $f$ satisfies
$\lambda\trianglerighteq\alpha$ and $\alpha_{1}\leq\lambda_{1}$ for all $i$.
Apply $\pi_{\lambda_{1}}$ to $f$ then by Lemma \ref{resagtb} $\beta
\trianglelefteq\left(  \lambda_{2},\lambda_{3},\ldots,\lambda_{N}\right)  $
for each $x^{\beta}$ appearing in $\pi_{\lambda_{1}}f$. For $i\geq2$
$\omega_{i}$ commutes with $\pi_{\lambda_{1}}$ and by Theorem \ref{xiPieq}
$\pi_{\lambda_{1}}\xi_{i}f$ =$\xi_{i}^{\prime}\pi_{\lambda_{1}}f$ . Thus
$\omega_{i}\pi_{\lambda_{1}}f=\xi_{i}^{\prime}\pi_{\lambda_{1}}f$ for
$i\geq2.$Also $\left(  \lambda_{2},\lambda_{3},\ldots,\lambda_{N}\right)
\in\mathbb{N}_{0}^{N-1,+}$ is $\vartriangleright$-maximal in $\pi_{\lambda
_{1}}f$.
\end{proof}

The definition of RSYT has been slghtly modified to allow filling with
$2,3,\ldots,N$. The isotype of $\pi_{\lambda_{1}}f$ is $\tau^{\prime}:=\left(
nd_{1}-1,nd_{2}-1,\ldots,nd_{L}-1,r_{L+1}-1\right)  .$

\begin{theorem}
\label{singFn}In the notation of Theorem \ref{singLB} if $d_{2}\geq2$ then
there is a singular polynomial for the parameter $\varpi$ in $n\left(
d_{1}+1\right)  -1$ variables with $\lambda=\left(  \left(  md_{1}\right)
^{n},0^{nd_{1}-1}\right)  $, of isotype $\left(  n,nd_{1}-1\right)  .$
\end{theorem}

\begin{proof}
Apply Proposition \ref{ResXiPi} repeatedly, and by hypothesis $nd_{2}%
-1\geq2n-1>n$. The remaining RSYT is
\[
Y^{\prime}=%
\begin{bmatrix}
N & N-1 & \ldots & \ldots & \ldots & N-nd_{1}+2\\
N-nd_{1}+1 & \ldots & N-nd_{1}-n+2 &  &  &
\end{bmatrix}
,
\]
and has the $t^{C}$-vector $\left[  t^{n-2},t^{n-1},\ldots,1,t^{-1}%
,t^{nd_{1}-2},t^{nd_{1}-3},\ldots,t,1\right]  .$
\end{proof}

\section{\label{conclud}Concluding Argument}

Re-index the variables by replacing $d_{1}\geq2$ (implied by $d_{2}\geq2$) by
$d$, $N$ by $N=nd-1+n$ and
\[
Y^{\prime\prime}=%
\begin{bmatrix}
nd-1+n & nd-2+n & \ldots & \ldots & \ldots & n+1\\
n & \ldots & 1 &  &  &
\end{bmatrix}
.
\]

\begin{proposition}
Suppose $\lambda=\left(  d^{n},0^{nd-1}\right)  $ and $\gamma\in\mathbb{N}%
_{0}^{K}$ for some $K\geq N$ satisfies $\left\vert \gamma\right\vert =nd$ and
$C_{i}:n\left(  \lambda_{i}-\gamma_{i}\right)  =r_{\gamma}\left(  i\right)
-i$ for $1\leq i\leq K$ (setting $\lambda_{i}=0$ for $i>N$) then
$\gamma=\lambda$ or $\gamma=\beta:=\left(  0^{n},1^{nd}\right)  .$
\end{proposition}

\begin{proof}
By condition $C_{n+1}$ we have $\left(  r_{\gamma}\left(  n+1\right)
-n-1\right)  =-n\gamma_{n+1}$ so that $\gamma_{n+1}=1-\frac{1}{n}\left(
r_{\gamma}\left(  n+1\right)  -1\right)  \leq1$ and thus $\gamma_{n+1}=1$ or
$\gamma_{n+1}=0$. If $\gamma_{n+1}=1$ then $r_{\gamma}\left(  n+1\right)  =1$,
which implies $\gamma_{i}=0$ for $1\leq i\leq n$ and $\gamma_{i}\leq1$ for
$i>n+1.$ If $j>n$ and $\gamma_{j}=0$ then by $C_{j}$ $r_{\gamma}\left(
j\right)  =j=\#\left\{  k<=j:\gamma_{k}\geq0\right\}  +\#\left\{
k>j:\gamma_{k}>0\right\}  $ so that $k>j$ implies $\gamma_{k}=0$. Since
$\left\vert \gamma\right\vert =\left\vert \lambda\right\vert =nd$ we see that
$\gamma_{n+1}=1$ implies $\gamma^{+}=\left(  1^{nd}\right)  $ and in fact
$\gamma_{i}=1$ for $n+1\leq i\leq n\left(  d+1\right)  $, since $\gamma_{j}=0$
and $\gamma_{j+1}=1$ is impossible for any $j>n$. If $1\leq j\leq n$ then
$r_{\gamma}\left(  j\right)  =nd+j$ and $n\left(  \lambda_{i}-\gamma
_{i}\right)  =nd=r_{\gamma}\left(  j\right)  -j$, thus satisfying $C_{j}$. The
other conditions $C_{i}$ are verified similarly. Thus $\gamma=\beta$.

If $\gamma_{n+1}=0$ then $r_{\gamma}\left(  n+1\right)  =n+1$ and $\ell\left(
\gamma\right)  =n$. Suppose $1\leq j\leq n$ then $C_{j}$ states $n\left(
\lambda_{j}-\gamma_{j}\right)  =r_{\gamma}\left(  j\right)  -j$ and the bounds
$1\leq j,r_{\gamma}\left(  j\right)  \leq n$ imply $\left\vert r_{\gamma
}\left(  j\right)  -j\right\vert \leq n-1$ and thus $\gamma_{j}=\lambda_{j}$.
\end{proof}

\begin{corollary}
Suppose $\lambda=\left(  \left(  md\right)  ^{n},0^{nd-1}\right)
\in\mathbb{N}_{0}^{N,+}$. The coeffients of $M_{\lambda}\left(  x\right)  $
have no poles at $\varpi$.
\end{corollary}

\begin{proof}
$M_{\lambda}\left(  x\right)  $ is a nonzero multiple of $x^{\lambda}%
+\sum_{\beta\vartriangleleft\lambda}A_{\lambda,\beta}x^{\beta}$. For each
$\beta\vartriangleleft\lambda$ there is at least one index $j_{\beta}$ such
that $\zeta_{\lambda}\left(  i_{\beta}\right)  \neq\zeta_{\beta}\left(
i_{\beta}\right)  $ at $\varpi$ or else $q^{\lambda_{i}-\beta_{i}}t^{r_{\beta
}\left(  i\right)  -i}=1$ for all $i\leq N$. In this case by Lemma \ref{qmtn}
$\left(  \lambda_{i}-\beta_{i}\right)  =ms_{i},r_{\beta}\left(  i\right)
-i=ns_{i}$ for some $s_{i}\in\mathbb{Z}$. Set $\lambda^{\prime}=\frac{1}%
{m}\lambda,\beta^{\prime}=\frac{1}{m}\beta$ then $n\left(  \lambda_{i}%
^{\prime}-\beta_{i}^{\prime}\right)  =r_{\beta}\left(  i\right)  -i$ for all
$i$ and by the Proposition $\beta^{\prime}=\lambda^{\prime}$ or $\beta
^{\prime}=\left(  0^{n}.1^{nd}\right)  $ but the latter is impossible because
$\left(  0^{n},1^{nd}\right)  \notin N_{0}^{N}$. Finally (this works because
there is a triangular expansion $x^{\lambda}=cM_{\lambda}+\sum\limits_{\beta
\vartriangleleft\lambda}A_{\beta,\lambda}^{\prime}M_{\beta}$ which holds for
generic $\left(  q,t\right)  $)
\[
M_{\lambda}\left(  x\right)  =c\prod\limits_{\beta\vartriangleleft\lambda
}\frac{\xi_{i_{\beta}}-\zeta_{\lambda}\left(  i_{\beta}\right)  }{\zeta
_{\beta}\left(  i_{\beta}\right)  -\zeta_{\lambda}\left(  i_{\beta}\right)
}x^{\lambda}.
\]
This shows that the poles of $M_{\lambda}$ are of the form $q^{a}t^{b}-1=0$
and $\varpi$ is not a pole.
\end{proof}

\begin{proposition}
Suppose $f$ is as in Theorem \ref{singFn} then $f\left(  x\right)
=cM_{\lambda}\left(  x\right)  $ at $\varpi$ for some constant $c\neq0$.
\end{proposition}

\begin{proof}
By matching coefficients of $x^{\lambda}$ find $c$ so that $\mathrm{coeff}%
\left(  x^{\lambda},f-cM_{\lambda}\right)  =0$. If $g:=f-cM_{\lambda}\neq0$
then there exists $\beta$ such that $x^{\beta}$ is $\vartriangleright$-maximal
in $g$. By $\vartriangleright$-triangularity $\xi_{i}g=q^{\beta_{i}%
}t^{N-r_{\beta}\left(  i\right)  }g$ (at $\varpi$) for all $i$. But $g$ has
the same eigenvalues as $M_{\lambda}$, that is, $q^{\beta_{i}}t^{N-r_{\beta
}\left(  i\right)  }=q^{\lambda_{i}}t^{N-i}$ at $\varpi$ and the proof of the
Corollary showed that $\beta=\lambda$ , contradicting $g\neq0$.
\end{proof}

Recall the transformation formula \ref{TMsi} for $M_{\alpha}$ for $\alpha
_{i}>\alpha_{i+1}$ with $z=\frac{\zeta_{\alpha}\left(  i+1\right)  }%
{\zeta_{\alpha}\left(  z\right)  }$%
\[
M_{s_{i}\alpha}=\frac{\left(  1-z\right)  ^{2}}{\left(  1-zt\right)  \left(
t-z\right)  }\left(  T_{i}+\frac{1-t}{1-z}\right)  M_{\alpha}.
\]
If $M_{\alpha}$ has no pole at $\varpi$ and $z\neq1,t,t^{-1}$ then
$M_{s_{i}\alpha}$ has no pole at $\varpi$. When $\alpha^{+}=\lambda$ then
$\alpha_{i}>\alpha_{i+1}$ implies $\alpha_{i}=md$ and $\alpha_{i+1}%
=0,z=q^{-md}t^{r_{\alpha}\left(  i\right)  -r_{\alpha}\left(  i+1\right)
}=t^{nd+r_{\alpha}\left(  i\right)  -r_{\alpha}\left(  i+1\right)  }$ at
$\varpi$. In the substring $\left(  \alpha_{1},\ldots,\alpha_{i},\alpha
_{i+1}\right)  $ there are $r_{\alpha}\left(  i\right)  $ values $md$ and
$i+1-r_{\alpha}\left(  i\right)  $ zeros, thus $r_{\alpha}\left(  i+1\right)
=n+i+1-r_{\alpha}\left(  i\right)  $. Thus $z=t^{b}$ with $b=nd+2r_{\alpha
}\left(  i\right)  -n-i-1$. Suppose $r_{\alpha}\left(  i\right)  =n$, thus
$i\geq n$ and $s_{i}$ can act on $\alpha$ without introducing a pole at
$\varpi$ if $nd+n-i-1>1$, that is $i<nd+n-2=N-1$. The last permitted
occurrence of $md$ in $\alpha$ is $i=N-2.$ Next move the second last
occurrence of $md$ in $\alpha$ as far as possible without a pole: set
$r_{\alpha}\left(  i\right)  =n-1$ and require $nd+2\left(  n-1\right)
-n-i-1>1$, that is, $i<nd+n-4=N-3$, thus $i=N-4$ is the last permitted value.
More generally let $r_{\alpha}\left(  i\right)  =n-j$ (with $0\leq j\leq n-1$)
then require $nd+2\left(  n-j\right)  -n-i-1>1$, that is, $nd+n-2j-2>i$ or
$i<N-1-2j$; the last permitted value is $i=N-2\left(  j+1\right)  .$

Let
\begin{align*}
\alpha &  =\left(  0^{nd-n-1},md,0,md,0.\ldots,md,0\right) \\
\zeta_{\alpha}  &  =\left[  t^{N-n-1},\ldots,t^{n},q^{md}t^{N-1}%
,t^{n-1},\ldots,q^{md}t^{N-n},1\right]  .
\end{align*}
We showed that $M_{\alpha}$ has no poles at $\varpi$, and if $M_{\lambda}$ at
$\varpi$ is singular then so is $M_{\alpha}$. The spectral vector
$\zeta_{\alpha}$ at $\varpi$ coincides with the $t^{C}$-vector of the RSYT%
\[
Y_{0}=%
\begin{bmatrix}
N & N-2 & \cdots & N-2n+2 & N-2n & \cdots & 1\\
N-1 & N-3 & \cdots & N-2n+1 &  &  &
\end{bmatrix}
,
\]
and thus $\omega_{N-1}Y_{0}=t^{-1}Y_{0}$; by construction $\zeta_{\alpha
}\left(  N-1\right)  =q^{md}t^{N-n}=t^{-nd+N-n}=t^{-1}$. If $M_{\alpha}$ at
$\varpi$ is singular then $\omega_{N-1}M_{\alpha}=\xi_{N-1}M_{\alpha}%
=t^{-1}M_{\alpha}$; this means%
\begin{align*}
t^{-1}T_{N-1}T_{N-1}M_{\alpha}  &  =t^{-1}M_{\alpha}\\
\left(  \left(  t-1\right)  T_{N-1}+t\right)  M_{\alpha}  &  =M_{\alpha}\\
\left(  t-1\right)  T_{N-1}M_{\alpha}  &  =\left(  1-t\right)  M_{\alpha}\\
\left(  T_{N-1}+1\right)  M_{\alpha}  &  =0.
\end{align*}
For the next step we recall some standard definitions: the $q$-Pochhammer
symbol is $\left(  a;q\right)  _{k}=\prod\limits_{i=1}^{k}\left(
1-aq^{i-1}\right)  $ and the generalized $\left(  q,t\right)  $-Pochhammer
symbol for $\lambda\in\mathbb{N}_{0}^{N,+}$ is%
\[
\left(  v;q,t\right)  =%
%TCIMACRO{\dprod _{i=1}^{N}}%
%BeginExpansion
{\displaystyle\prod_{i=1}^{N}}
%EndExpansion
\left(  vt^{1-i};q\right)  _{\lambda_{i}}.
\]
In the context of the Ferrers diagram representation of a composition
$\alpha\in\mathbb{N}_{0}^{N}$, $\left\{  \left(  i,j\right)  :1\leq i\leq
N,1\leq j\leq\alpha_{i}\right\}  $ (the rows with $\alpha_{i}=0$ are empty)
define the arm-length and leg-length of a box in the diagram ($\lambda
\in\mathbb{N}_{0}^{N,+}$)
\begin{align*}
\mathrm{arm}\left(  i,j;\lambda\right)   &  :=\lambda_{i}-j,\\
\mathrm{arm}\left(  i,j;\alpha\right)   &  :=\alpha_{i}-j,\\
\mathrm{leg}\left(  i,j;\lambda\right)   &  :=\#\left\{  l:i<l\leq
N,j\leq\lambda_{l}\right\}  ,
\end{align*}

\[
\mathrm{leg}\left(  i,j;\alpha\right)  :=\#\left\{  r:r>i,j\leq\alpha_{r}%
\leq\alpha_{i}\right\}  +\#\left\{  r:r<i,j\leq\alpha_{r}+1\leq\alpha
_{i}\right\}  .
\]
The $\left(  q,t\right)  $-hook product is
\[
h_{q,t}\left(  v;\alpha\right)  =\prod\limits_{\left(  i,j\right)  \in\lambda
}\left(  1-vq^{\mathrm{arm}\left(  i,j;\alpha\right)  }t^{l\mathrm{eg}\left(
i,j;\alpha\right)  }\right)  .
\]
There is an evaluation at a special point (see \cite[Cor. 7]{DL2015}): let
$x^{\left(  0\right)  }:=\left(  1,t,t^{2},\ldots,t^{N-1}\right)  $, then for
any $\beta\in\mathbb{N}_{0}^{N}$%
\[
M_{\beta}\left(  x^{\left(  0\right)  }\right)  =q^{b\left(  \beta\right)
}t^{e^{\prime}\left(  \beta^{+}\right)  }\frac{\left(  qt^{N};q,t\right)
_{\beta^{+}}}{h_{q,t}\left(  qt;\beta\right)  },
\]
where $b\left(  \beta\right)  =\sum_{i=1}^{N}\binom{\beta_{i}}{2},e^{\prime
}\left(  \beta^{+}\right)  =\sum_{i=1}^{N}\beta_{i}^{+}\left(  N-i\right)  $.

\begin{theorem}
$\left(  T_{N-1}+1\right)  M_{\alpha}\neq0$ at $\varpi$ and $M_{\alpha}$ is
not singular.
\end{theorem}

\begin{proof}
For any polynomial $p$ let $x=x^{\left(  0\right)  }$ in $T_{i}p\left(
x\right)  =\left(  1-t\right)  x_{i+1}\frac{p\left(  x\right)  -p\left(
xs_{i}\right)  }{x_{i}-x_{i+1}}+tp\left(  xs_{i}\right)  $ then $T_{i}p\left(
x^{\left(  0\right)  }\right)  =t\left(  p\left(  x^{\left(  0\right)
}\right)  -p\left(  x^{\left(  0\right)  }s_{i}\right)  \right)  +tp\left(
x^{\left(  0\right)  }s_{i}\right)  =tp\left(  x^{\left(  0\right)  }\right)
$ (since $x_{i+1}^{\left(  0\right)  }=tx_{i}^{\left(  0\right)  }$). Set
$b_{0}=b\left(  \alpha\right)  =n\binom{md}{2}$, $e_{0}=e^{\prime}\left(
\alpha^{+}\right)  =\frac{1}{2}mdn\left(  2N-n-1\right)  $ then%
\begin{gather*}
T_{N-1}M_{\alpha}\left(  x^{\left(  0\right)  }\right)  +M_{\alpha}\left(
x^{\left(  0\right)  }\right)  =\left(  t+1\right)  M_{\alpha}\left(
x^{\left(  0\right)  }\right) \\
=q^{b_{0}}t^{e_{0}}\left(  t+1\right)  \frac{\left(  q^{N}t;q,t\right)
_{\alpha^{+}}}{h_{q,t}\left(  qt;\alpha\right)  }.
\end{gather*}
The numerator is%
\[
\left(  q^{N}t;q,t\right)  _{\alpha^{+}}=\prod\limits_{i=1}^{n}\left(
qt^{N-i+1};q\right)  _{md}=\prod\limits_{i=1}^{n}\prod_{j=1}^{dm}\left(
1-q^{j}t^{nd+n-i}\right)  ,
\]
where the only term vanishing at $\varpi$ is for $i=n,j=dm$ (for suppose
$j=rm$ with $r\leq d,nd+n-i=rn$ for some $r\in\mathbb{N}$ then $n\geq
i=n\left(  d-r+1\right)  $ and $d-r+1\leq1$, that is, $r\geq d$, hence
$r=d,i=n$). For the hook product observe that if $1\leq j\leq n$ then
$\mathrm{leg}\left(  \alpha;N-2j+1,1\right)  =nd-2$ because there are $nd-1-j$
zero values in $\left(  \alpha_{1},\ldots,\alpha_{N-2j+1}\right)  $ and $j-1$
values of $md$ in $\left(  \alpha_{N-2j+2},\ldots,\alpha_{N}\right)  $. Since
$\mathrm{arm}\left(  \alpha;N-2j+1,1\right)  =dm-1$ we find that the boxes
$\left\{  \left[  N-2j+1,1\right]  :1\leq j\leq n\right\}  $ contribute
$\left(  1-q^{dm}t^{nd-1}\right)  ^{n}$ to $h_{q,t}\left(  qt;\alpha\right)
$. This term becomes $\left(  1-t^{-1}\right)  ^{n}$ at $\varpi$. The other
boxes in the diagram of $\alpha$ are $\left\{  \left[  N-2j+1,k\right]  :1\leq
j\leq n,2\leq k\leq md\right\}  $ and $\mathrm{leg}\left(  \alpha
;N-2j+1,k\right)  =j-1$, $\mathrm{arm}\left(  \alpha;N-2j+1,k\right)  =dm-k$.
Thus%
\begin{align*}
h_{q,t}\left(  qt;\alpha\right)   &  =\left(  1-q^{dm}t^{nd-1}\right)
^{n}\prod\limits_{j=1}^{n}\prod\limits_{k=1}^{dm}\left(  1-q^{dm-k+1}%
t^{j}\right) \\
&  =\left(  1-q^{dm}t^{nd-1}\right)  ^{n}\prod\limits_{j=1}^{n}\prod
\limits_{i=1}^{dm}\left(  1-q^{i}t^{j}\right)  .
\end{align*}
The only term in the product vanishing at $\varpi$ is for $i=m,j=n$. Thus the
term $\left(  1-q^{m}t^{n}\right)  $ cancels out in $\frac{\left(
q^{N}t;q,t\right)  _{\alpha^{+}}}{h_{q,t}\left(  qt;\alpha\right)  }$ and
$\left(  T_{N-1}+1\right)  M_{\alpha}\left(  x^{\left(  0\right)  }\right)
\neq0$.
\end{proof}

\begin{example}
Let $N=5$, $n=2$, $m=1$, $d=2$ then $\alpha=\left(  0,2,0,2,0\right)  $ and
$\varpi=\left(  t^{-2},t\right)  $ (that is, $qt^{2}=1$) The spectral vector
of $\alpha$ is $\left[  t^{2},q^{2}t^{4},t,q^{2}t^{3},1\right]  $ which equals
$\left[  t^{2},1,t,t^{-1},1\right]  $ at $q=t^{-2}$. The expression for
$M_{\alpha}$ is too large to display here (32 monomials); the denominators of
the coefficients are factors of $qt-1,\left(  q^{2}t^{3}-1\right)  ^{2}$ and
\[
M_{\alpha}\left(  1,t,t^{2},t^{3},t^{4}\right)  =q^{2}t^{14}\frac{\left(
qt^{2}+1\right)  \left(  qt^{4}-1\right)  \left(  qt^{5}-1\right)  \left(
q^{2}t^{5}-1\right)  }{\left(  q^{2}t^{3}-1\right)  ^{2}\left(  qt-1\right)  }%
\]
which does not vanish at $q=t^{-2}$. But the same monomial is singular with
$n=4$, $d=1$, $m=2$ and $q=-t^{-2}$ (that is, $q^{2}t^{4}=1$ but $qt^{2}\neq
1$). The singularity can be proven by direct computation and the vanishing of
$M_{\alpha}\left(  1,\ldots,t^{4}\right)  $ is only a necessary condition.
\end{example}

We have shown if there is a singular polynomial as described in Theorem
\ref{singLB} and $d_{2}\geq2$ then by using the restriction Proposition
\ref{ResXiPi} repeatedly there is a singular polynomial of isotype $\left(
nd_{1}-1,n\right)  $, which in turn implies that $M_{\alpha}$ is singular.
This is impossible and we conclude that $d_{2}=1$ is necessary, and all
singular polynomials have been determined.

\end{document}